\documentclass[11pt,reqno,oneside]{amsart}
\input epsf

\def\Box{\setlength{\unitlength}{0.01cm}
     \begin{picture}(15,15)(-5,-5)
       \framebox(15,15){}
     \end{picture} }

\newtheorem{theorem}{Theorem}[section]
\newtheorem{lemma}{Lemma}[section]

\newtheorem{remark}{Remark}[section]

\def\proof{\noindent{\bf Proof.}\,\,}

\def\O{\Omega}
\def\R{{\mathbb R}}

\def\d{\delta}
\def\g{\gamma}

\begin{document}
\section*{}

\setcounter{equation}{0}
\title[]{Improved Poincar\'e inequalities with weights}

\author{Irene Drelichman}
\address{Departamento de Matem\'atica, Facultad de Ciencias Exactas y Naturales,
Universidad de Buenos Aires, 1428 Buenos Aires, Argentina} \email{irene@drelichman.com}

\author{Ricardo G. Dur\'an}
\address{Departamento de Matem\'atica, Facultad de Ciencias Exactas y Naturales,
Universidad de Buenos Aires, 1428 Buenos Aires, Argentina} \email{rduran@dm.uba.ar}

\thanks{Supported by ANPCyT under grant  PICT 03-13719, by Universidad de Buenos Aires under grant X052 and by
CONICET under grant PIP 5478. The authors are members of CONICET,
Argentina.}

\keywords{weighted Sobolev inequality, weighted Poincar\'e inequality, reverse doubling weights, John domains}
\subjclass[2000]{46E35, 26D10}

\begin{abstract}
In this paper we prove that if $\O\in\mathbb{R}^n$ is a bounded John domain, the
following weighted Poincar\'e-type inequality holds:
$$
\inf_{a\in \mathbb{R}}\| (f(x)-a) w_1(x) \|_{L^q(\Omega)} \le C
\|\nabla f(x) d(x)^\alpha w_2(x) \|_{L^p(\O)}
$$
where $f$ is a locally Lipschitz function on $\O$, $d(x)$ denotes the distance of $x$ to the boundary of $\O$,
the weights $w_1, w_2$ satisfy certain cube conditions, and $\alpha \in [0,1]$ depends on $p,q$ and $n$.
This result generalizes previously known weighted inequalities, which can also be obtained with
our approach.

\end{abstract}
\maketitle

\section{Introduction}

The purpose of this paper is to present a simple unified approach
to prove weighted Poincar\'e-type inequalities in John domains.

The class of John domains was first introduced in \cite{J} and
named after the author of that paper by Martio and Sarvas
\cite{MS}. It contains Lipschitz domains as well as other domains
with very non-regular boundaries, and it has played an important role
in several problems in analysis. In particular, as it has been made
clear in \cite{BK}, it is closely connected to the improved
Poincar\'e inequalities we are interested in.

The Sobolev-Poincar\'e inequality
\begin{equation}
\label{sob-poinc} \inf_{a \in \mathbb{R}} \| f(x) -
a\|_{L^{\frac{np}{n-p}}(\Omega)} \le C \|\nabla f(x) \|_{L^p(\Omega)}
\end{equation}
with $\Omega \subseteq \mathbb{R}^n$ being a John domain, and $f\in W^{1,p}(\Omega)$, was  proved in the case $1<p<n$ in
\cite{M}, and later extended to the case $p=1$ in \cite{Bo}. See
also \cite{H} for proofs, other references and a nice
account on the history of this problem.

Moreover, it was proved in \cite{BK} that John domains are
essentially the largest class of domains for which this inequality
can hold, more precisely, if $\Omega \subseteq \mathbb{R}^n$ is a
domain of finite volume that satisfies a separation property
(cf.\cite{BK}) and $1\le p<n$, then $\Omega$ satisfies the
Sobolev-Poincar\'e inequality if and only if it is a John
domain.

The Sobolev-Poincar\'e inequality can be seen as a special case of
a much wider family of so-called improved Poincar\'e inequalities.
Indeed, it was proved in \cite{R} that if $\Omega\subseteq
\mathbb{R}^n$ is a bounded John domain, and $f\in L_{loc}^1
(\Omega)$ is such that $\nabla f(x) d(x)^\alpha \in L^p(\Omega)$,
then

\begin{equation}
\label{improv-poinc}
 \inf_{a \in \mathbb{R}} \| f(x) - a\|_{L^q(\Omega)} \le C
\|\nabla f(x) d(x)^\alpha\|_{L^p(\Omega)}
\end{equation}
whenever
$1<p\le q\le \frac{np}{n-p(1-\alpha)}$ with $p(1-\alpha)<n$, and
$\alpha \in [0,1]$, with $d(x)$ being the distance of a point $x$
to the boundary of $\Omega$ (the same inequality holds for
unbounded John domains with $1\le p\le q = \frac{np}{n-p(1-\alpha)}$).
Letting $\alpha=0$ in (\ref{improv-poinc}) one clearly obtains
inequality (\ref{sob-poinc}).

A further generalization of Poincar\'e inequalities in weighted
spaces was made in \cite{C} for bounded John domains. It was shown in that paper
that under certain cube conditions on the weights $w_1, w_2$, the
following inequality holds for bounded John domains:

\begin{equation}
\label{poinc-pesos}
 \inf_{a \in \mathbb{R}} \| (f(x) - a) w_1(x)\|_{L^q(\Omega)} \le C
\|\nabla f(x) w_2(x)\|_{L^p(\Omega)}
\end{equation}
whenever $f$ is a Lipschitz function and $1<p\le q<\infty$. Notice
that the author of \cite{C} refers to domains satisfying the Boman
chain condition, but for connected domains in $\mathbb{R}^n$ this
is exactly the same class as that of John (see \cite{B} for proof
of this inequality even in a much more general context).

Inequality (\ref{poinc-pesos}) can also be extended to unbounded John domains as
it was done in \cite{R} for the case of (\ref{improv-poinc}) (see \cite{R2}). Both
results rely heavily on the main theorem of \cite{V}, which states that an unbounded John
domain can be written as an increasing union of bounded John domains in a way that
allows to pass to the limit using the dominated convergence theorem.

As we did for inequality (\ref{sob-poinc}), we could also think of
 inequality (\ref{poinc-pesos}) as a special case of  a wider
family of inequalities explicitly involving powers of the distance
to the boundary. Indeed, we will prove in this paper that if $f$ is a locally
Lipschitz function on $\Omega$

\begin{equation}
\label{improv-pesos}
 \inf_{a \in \mathbb{R}} \| (f(x) - a) w_1(x)\|_{L^q(\Omega)} \le C
\|\nabla f(x) d(x)^\alpha w_2(x) \|_{L^p(\Omega)}
\end{equation}
for suitable weights $w_1, w_2$, and with $\alpha$ depending on
$p, q$ as in inequality (\ref{improv-poinc}), thus extending the
results in \cite{C}. Notice that when the density of the locally Lipschitz functions
 in the involved weighted norms holds, this result extends to functions
in the corresponding weighted Sobolev spaces.

It is worth noting that the technique we will use for the proof of
 inequality (\ref{improv-pesos}) differs completely from the one used  in \cite{C} for the case $\alpha=0$.
Instead of relying on chains of cubes and cube-by-cube
inequalities, we recover the simpler classical ideas which relate
Sobolev-Poincar\'e inequalities with fractional integrals (see,
e.g., \cite{H} and references therein). Similar ideas were
previously used for John domains in \cite{M} to prove the
Sobolev-Poincar\'e inequality, but the fact that they can also be
used in connection with the distance to the boundary seems to be
new.

We will use a representation formula proved in \cite{ADM} that
essentially allows us to recover $f$ from its gradient (an
alternative proof of inequality (\ref{sob-poinc}) can also be
found in that paper). It has, as mentioned before, the advantage
of allowing us to introduce the distance to the boundary without
recurring to Whitney cubes, and it will allow us to reduce the
proof of inequalities (\ref{improv-poinc}) and
(\ref{improv-pesos}) to known continuity results for fractional
integrals and the Hardy-Littlewood maximal function.

Although inequality (\ref{improv-poinc}) can be seen as a special
case of (\ref{improv-pesos}) taking $w_1=w_2=1$, we have chosen to
present them separately for the sake of clarity and because the
hypotheses needed are weaker than those we require for the more general
cases. We will also split inequality (\ref{improv-pesos}) into the
cases $w_1 = w_2$ and $w_1 \not= w_2$. We shall refer to the first
case as `one-weighted' case and to
the second one as `two-weighted' case. Once the ideas are made
clear in the simpler cases, we shall be somewhat sketchy to
indicate how they can be adapted to the more general case.

The paper is organized as follows. In section 2 we recall some
definitions, obtain the representation formula that we will be
using in the remainder of the paper, and show how it relates to
the distance to the boundary. Section 3 is devoted to the
unweighted and one-weighted cases. We obtain a simpler proof of
the results in \cite{R} and, following the technique presented in
\cite{H}, we extend inequality (\ref{improv-pesos}) for $w_1=w_2$
to the previously unknown case $p=1$ (Theorem \ref{teo34}).
Finally, in section 4 we show how our arguments can be used to generalize the results in \cite{C} and obtain new
inequalities in the two-weighted case (Theorems \ref{teo41} and \ref{teo42}).

\section{Preliminaries}

The notation used in this paper is rather standard. By $C$ we will denote a general constant which can change its value even within a single string of estimates. We will write $C(*,...,*)$ to emphasize that the constant depends on the quantities appearing in the parentheses only.

By a weight function we mean a nonnegative measurable function on $\R^n$.

 Given a domain $\O\subset\R^n$ and any
$x\in\O$, we let $d(x)$ denote the distance of $x$ to the boundary of
$\O$. A bounded domain $\O\subset\R^n$ is a John domain if for a
fixed $x_0\in\O$ and any $y\in\O$ there exists a rectifiable
curve, called John curve, given by
$$
\g(\cdot,y)\,:\, [0,1]\to\O
$$
such that $\g(0,y)=y$ and $\g(1,y)=x_0$, and there exist constants
$\d$ and $K$, depending only on the domain $\O$ and on $x_0$, such
that
\begin{equation}
\label{john1} d(\g(s,y))\ge\d s
\end{equation}
and
\begin{equation}
\label{john2} |\dot\g(s,y)|\le K
\end{equation}
where $\dot\g(s,y):=\frac{\partial\g}{\partial s}(s,y)$.

 In what follows, we will be using that $\g(s,y)$ and
$\dot\g(s,y)$ are measurable functions. This property need not be
fulfilled if we take $\g(\cdot, y)$ to be an arbitrary John curve
for each fixed $y\in\Omega$, but it can be obtained by means of a
slight technical modification of a given family of curves (see \cite[Lemma 2.1]{ADM} for
details). Moreover, to simplify notation we will assume,
without loss of generality, that $x_0=0$.

Let $\varphi\in C_0^\infty$ such that $\int_{\O}\varphi =1$ and
$\mbox{supp\,}\varphi \subset B(0,\delta/2)$. Given a locally Lipschitz function
$f$, we denote by $f_\varphi$ the weighted average of
$f$, namely, $f_\varphi=\int_\O f\varphi$.

The following lemmas of this section will be fundamental for the remainder of
this paper. They were proved in \cite{ADM} but
we have chosen to reproduce their proofs here for the sake of completeness.

\begin{lemma}
With the above notations, if $\O\subseteq \mathbb{R}^n$ is a John
domain and $y\in \O$,
\begin{equation}
\label{green}
f(y)-f_\varphi= \int_\O G(x,y)\cdot\nabla f(x)\, dx
\end{equation}
with
\begin{equation}
\label{Gj}
G(x,y):=-\int_0^1\left(\dot\gamma(s,y)+\frac{x-\gamma(s,y)}{s}\right)
\varphi\left (\frac{x-\gamma(s,y)}{s}\right)\frac{1}{s^n}\, ds.
\end{equation}
\end{lemma}

\proof
In view of (\ref{john1}), for any $y\in\O$ and $z\in
B(0,\delta/2)$ the curve given by
$$
\g(s,y) + sz \quad , \quad s\in[0,1]
$$
which joins $y$ and $z$, is contained in $\O$. Then
$$
f(y)-f(z)= - \int_0^1 \nabla f(\g(s,y)+sz)\cdot (\dot\g(s,y)+z)\,
ds.
$$
Multiplying by $\varphi(z)$ and integrating in $z$ we obtain
$$
f(y)- f_\varphi= - \int_\O \int_0^1 \nabla f(\g(s,y)+sz)\cdot
(\dot\g(s,y)+z)\varphi(z)\, ds dz.
$$
Making the change of variable $x=\g(s,y)+sz$ we have
$$
f(y)-f_\varphi= - \int_0^1 \int_\O \nabla f(x) \cdot
\left(\dot\g(s,y)+ \frac{x-\g(s,y)}{s}\right)
\varphi\left(\frac{x-\g(s,y)}{s}\right)\frac{1}{s^n}\, dx ds
$$
as we wanted to prove. \Box

\begin{lemma}
\label{cotadeG} There exists a constant $C=C(n,\d,K)$ such that
\begin{equation}
\label{boundofG} |G(x,y)|\le
C\frac{\|\varphi\|_{\infty}}{|x-y|^{n-1}}.
\end{equation}
\end{lemma}
\proof If $(x-\gamma(s,y))/s\in \mbox{supp\,}\varphi$ then
$|x-\gamma(s,y)| < (\d/2) s$. Therefore, using (\ref{john2}) and
$\g(0,y)=y$ we have
\begin{equation}
\label{xmenosy} |x-y| \le |x-\gamma(s,y)| +
|\gamma(s,y)-\gamma(0,y)| \le (\d/2) s + K s.
\end{equation}

Therefore,

$$
G(x,y) = \int_{C |x-y|}^1 \left\{
\dot{\gamma}(s,y)+\frac{x-\gamma(s,y)}{s} \right \} \varphi \left
(\frac{x-\gamma(s,y)}{s}\right)\frac{1}{s^n}\, ds .
$$
And, since

$$
\left|\dot{\gamma}(s,y)+\frac{x-\gamma(s,y)}{s}\right| \le K +
{\d/2},$$ the above estimate follows easily. \Box

\begin{lemma}
\label{cotax-y} There exists a constant $C=C(\d,K)$ such that, if
$G(x,y) \not= 0$, then
$|x-y|\le C(\d,K) d(x)$.
\end{lemma}

\proof Notice that, if $G(x,y)\not= 0$, there exists $s$ such that $\varphi\left(\frac{x-\g(s,y)}{s}\right) \not = 0$.

Let $\bar x \in \partial\Omega$ be
such that $d(x)=|x- \bar x|$. By (\ref{xmenosy}) and property (\ref{john1}), we have

$$
|x-y| \le \left(\frac{\d}{2} + K\right) s \le \left(\frac12
+\frac{K}{\d}\right) \d s \le   \left(\frac12 +\frac{K}{\d}\right)
d(\gamma(s,y)).
$$

But,

$$
d(\gamma(s,y))\le |\gamma(s,y)- \bar x|\le |\gamma(s,y) - x| + |x
- \bar x|\le \frac{\d s}{2} + d(x) \le \frac{d(\gamma(s,y))}{2} +
d(x),
$$

whence

$$
|x-y| \le  \left(1 +\frac{2K}{\d}\right) d(x).
$$
\Box

\section{The unweighted and one-weighted cases}

Since the case $p=1$ of the inequalities we are considering is different in nature from the remaining values
of $p$, we will split the proof of both the weighted and unweighted cases into two theorems, respectively.

\begin{theorem}
If $\O\subseteq \mathbb{R}^n$ is a bounded John domain,
\begin{equation}
\label{teo1}
 \inf_{a\in \mathbb{R}} \| f-a\|_{L^q(\O)} \le C \|
\nabla f(x) d(x)^\alpha\|_{L^p(\Omega)}
\end{equation}
whenever $f\in L^q(\O)$ is a locally Lipschitz function, $\alpha \in [0,1]$, $p(1-\alpha) <n$, and $1<p\le q\le
\frac{np}{n-p(1-\alpha)}$.
\end{theorem}

\proof
By duality,

$$ \| f - f_\varphi\|_{L^q(\O)} = \sup_{g\in
L^{q'}(\O)} \frac{\int_\O (f-f_\varphi) g}{\| g\|_{L^{q'}(\O)}},
$$
with $q'$ being the dual
exponent of $q$, $1/q + 1/q'=1$. Therefore, it suffices to obtain a bound for $\int_\O (f
-f_\varphi)g$ for $g\in L^{q'}(\O)$.

Using the representation formula (\ref{green}), we can write

$$
\int_\O (f(y) - f_\varphi) g(y) \, dy = \int_\O \int_\O
G(x,y)\cdot\nabla f(x)\, dx \, g(y)\, dy
$$

Interchanging the order of integration and using lemmas
\ref{cotadeG} and \ref{cotax-y}, we obtain

\begin{equation}
\label{dualidad}
 \int_\O |(f(y) - f_\varphi) g(y)| \, dy \le C
\int_\O \int_{|x-y|\le C d(x)} \frac{|g(y)|}{|x-y|^{n-1}} \, dy
|\nabla f(x)| \, dx
\end{equation}

We consider separately the cases $\alpha \in [0,1)$ and
$\alpha=1$.

In the case $\alpha \in [0,1)$, if we denote $I_\beta g(x)= \int
g(y) |x-y|^{\beta-n}\, dy$, we can bound the above expression by

\begin{equation}
\label{intfrac}
 C \int_{\mathbb{R}^n}|\nabla f (x)| d(x)^\alpha
I_{1-\alpha}|g(x)| \, dx
\end{equation}
where we have assumed that $|\nabla f|$ and $g$ are extended by
zero outside $\O$. Applying H\"older's inequality and the
continuity of the fractional integral (see, e.g., \cite{S}), this
expression can be bounded by

$$
 \| \nabla f(x) d(x)^\alpha\|_{L^p(\O)} \| I_{1-\alpha}
|g(x)|\|_{L^{p'}(\mathbb{R}^n)} \le C \| \nabla f(x)
d(x)^\alpha\|_{L^p(\O)} \| g(x)\|_{L^{q'}(\O)}
$$
thus proving (\ref{teo1}) in the case $\alpha \in [0,1)$.

In the case $\alpha=1$ (that is, $p=q$), a standard calculation
(see, e.g., \cite[Lemma 2.8.3]{Z}) shows that (\ref{dualidad}) can
be bounded by
$$
C \int_\Omega Mg(x) d(x) |\nabla f(x)| \, dx \le C
\|Mg(x)\|_{L^{q'}(\O)} \| \nabla f(x) d(x)\|_{L^q(\O)},
$$
and the desired result follows by boundedness of the
Hardy-Littlewood maximal function in $L^{q'}(\O)$ (see, e.g.,
\cite{S}). \Box

\begin{theorem}
If $\O\subseteq \mathbb{R}^n$ is a bounded John domain,
\begin{equation}
\label{teo1-L1}
 \inf_{a\in \mathbb{R}} \| f-a\|_{L^{n/(n-1+\alpha)}(\O)} \le C \|
\nabla f(x) d(x)^\alpha\|_{L^1(\Omega)}
\end{equation}
whenever $f\in L^{n/(n-1+\alpha)}(\O)$ is a locally Lipschitz function,  $1-\alpha<n$, and $\alpha \in [0,1]$.
\end{theorem}

\proof
In the case $\alpha=1$, inequality (\ref{teo1-L1}) can be proved
as in the previous theorem, using the continuity of the maximal
function in $L^\infty(\Omega)$.

In the case $\alpha\in [0,1)$, we follow the approach used in
\cite{H} to prove the Sobolev-Poincar\'e inequality for John
domains, modifying it to include the distance to the boundary in
our estimates.

For $g\in L^1(\O)$, let

$$
E_t=\left\{ x \in \O: \int_{\O} \frac{g(y)}{|x-y|^{n-1+\alpha}} \,dy >
t\right\}
$$

Then,
$$
|E_t| \le \int_E \int_\O \frac{g(y)}{t |x-y|^{n-1+\alpha}} \, dy \,
dx
$$

But,

$$
\int_{E_t} \frac{1}{ |x-y|^{n-1+\alpha}} \, dx \le C |E_t|^{(1-\alpha)/n}
$$
(see, e.g.,  \cite[inequality 7.2.6]{Jo}). Therefore,

$$
|E_t| t^{n/(n-1+\alpha)} \le C \left( \int_\O |g(y)| \, dy
\right)^{n/(n-1+\alpha)}
$$

Since, as in the proof of (\ref{teo1}),

$$
|f- f_\varphi| \le C \int_\O \frac{|\nabla f(y)|
d(y)^\alpha}{|x-y|^{n-1+\alpha}} \, dy,
$$
we conclude that

$$
\sup_{t>0} \left| \left\{ x\in \O : |f-f_\varphi| > t
\right\}\right| t^{n/(n-1+\alpha)} \le C \left( \int_\O |\nabla
f(y)| d(y)^\alpha \, dy \right)^{n/(n-\alpha)}
$$

This in turn implies, by \cite[Theorem 4]{H}, that

$$
\inf_{a\in \mathbb{R}} \| f(x) - a \|_{L^{n/(n-1+\alpha)}(\O)} \le
C \| \nabla f(x) d(x)^\alpha\|_{L^1(\O)}
$$

\Box

\begin{theorem}
Let $\O\subseteq\mathbb{R}^n$ be a bounded John domain. If $w$ is
a nonnegative function such that there exists a constant
$K<\infty$ such that

\begin{equation}
\label{cond-1peso}
\left( \frac{1}{|Q|}\int_Q w(x)^q \, dx
\right)^{1/q} \left( \frac{1}{|Q|} \int_Q w(x)^{-p'} \, dx
\right)^{1/p'} \le K
\end{equation}
where $Q$ is any $n$ dimensional cube, and $K$ is independent of
$Q$, then

$$
\inf_{a\in \mathbb{R}}\| (f(x)-a) w(x)\|_{L^q(\Omega)} \le C
\|\nabla f(x) d(x)^\alpha w(x)\|_{L^p(\O)}
$$
for all locally Lipschitz $f$, where $0< \alpha \le
1$, $p(1-\alpha)< n$ and $1< p \le q \le \frac{np}{n-p(1-\alpha)}$.
\end{theorem}

\proof
By duality, it suffices to bound $\int_\O (f-f_\varphi)(y) g(y)
\, dy$ for any $g$ such that $\|g(x) w^{-1}(x)\|_{L^{q'}(\O)}<
\infty$.

In the case $\alpha \in [0,1)$, using, as before, the bound (\ref{intfrac})
and H\"older's inequality, we obtain

$$
 \int_\O |(f- f_\varphi)(y) g(y)| \, dy \le  C
\| \nabla f(x) d(x)^\alpha w(x) \|_{L^p(\O)} \| I_{1-\alpha} |g(x)|
w(x)^{-1}\|_{L^{p'}(\O)}
$$

But, by condition (\ref{cond-1peso}), \cite[Theorem 4]{MW} and the fact that $I_{1-\alpha}$ is self-adjoint,

$$\|I_{1-\alpha} |g(x)| w^{-1}(x)\|_{L^{p'}(\O)} \le C \| g(x)
w^{-1}(x)\|_{L^{q'}(\O)}
$$
and the theorem follows.

In the case $\alpha=1$, bound (\ref{intfrac}), as before, by

$$
C \int_\O Mg(x) d(x) |\nabla f(x)| \, dx \le C \| Mg(x) w^{-1}(x)\|_{L^{p'}(\O)} \| \nabla f(x) d(x) w(x)\|_{L^p(\O)}
$$
and the result follows, since by condition (\ref{cond-1peso}) and \cite[Theorem 1.2]{CU}
(see also references therein for previously known results),
$$
\| Mg(x) w^{-1}(x)\|_{L^{p'}(\O)} \le C \| g(x) w^{-1}(x)\|_{L^{q'}(\O)}.
$$

\Box

\begin{remark}
Notice that if $w$ satisfies condition (\ref{cond-1peso}),  then $w^q$ belongs to Muckenhoupt's class
$A_r$ with $r=\frac{q}{p'}+1$, and therefore it is a doubling weight (which in turn implies that it satisfies the weaker `reverse doubling condition'
required for \cite[Theorem 1.2]{CU}).
\end{remark}

\begin{theorem}
\label{teo34}
Let $\O\subseteq\mathbb{R}^n$ be a bounded John domain. If $w$ is
a nonnegative function such that there exists a constant
$K<\infty$ such that

\begin{equation}
\label{cond-1peso-L1}
 \left( \frac{1}{|Q|} \int_Q w(x)^{\frac{n}{n-1+\alpha}} \, dx
\right)^{\frac{n-1+\alpha}{n}} \left( \mbox{ess }\sup_{x\in Q} \frac{1}{w(x)}
\right) < K
\end{equation}
where $Q$ is any $n$ dimensional cube and $K$ is independent of $Q$, then

$$
\inf_{a\in \mathbb{R}}\| (f(x)-a)
w(x)\|_{L^{n/(n-1+\alpha)}(\Omega)} \le C \|\nabla f(x)
d(x)^\alpha w(x)\|_{L^1(\O)}
$$
for all locally Lipschitz $f$ and $
\alpha \in [0,1)$.

When $\alpha=1$, condition (\ref{cond-1peso-L1}) should be replaced by

$$
M(w(x))\le C w(x)
$$
for almost every $x\in\O$ (that is, $w\in A_1$).
\end{theorem}

\proof
In the case $\alpha \in [0,1)$, for each $t>0$ let $E_t = \{ |I_{1-\alpha}g(x)|>t\}$. By \cite[Theorem 5]{MW}, if
$w$ satifies condition (\ref{cond-1peso-L1}),

$$
\int_{E_t} w(x)^{\frac{n}{n-1+\alpha}} \, dx \le C t^{-\frac{n}{n-1+\alpha}} \left( \int_{\mathbb{R}^n}
|g(x)| w(x) \, dx\right)^{\frac{n}{n-1+\alpha}}
$$

But, as before,

$$
|f- f_\varphi| \le C \int_\O \frac{|\nabla f(y)|
d(y)^\alpha}{|x-y|^{n-1+\alpha}} \, dy = C I_{1-\alpha} (|\nabla
f| d(x)^\alpha)
$$

Therefore, setting $d\mu = w(x)^{n/(n-1+\alpha)} \, dx$, we obtain
that

$$
\mu\{ |f-f_\varphi| > t\} t^{n/(n-1+\alpha)} \le C \mu \{
I_{1-\alpha} (|\nabla f| d(x)^\alpha)> t\} t^{n/(n-1+\alpha)}
$$

$$
\le C \left(\int_\O |\nabla f(x)| d(x)^\alpha w(x) \,
dx\right)^{n/(n-1+\alpha)}
$$
which, by \cite[Lemma 4]{H}, implies

$$
\inf_{a\in \mathbb{R}} \left( \int_\O |f-a|^{n/(n-1+\alpha)} \,
d\mu\right)^{(n-1+\alpha)/n} \le C \int_\O |\nabla f| \, d\nu
$$
where $d\nu = d(x)^\alpha w(x) \, dx$, that is,

$$
\inf_{a\in \mathbb{R}} \|(f-a)(x) w(x)\|_{L^{n/(n-1+\alpha)}(\O)}
\le C \| \nabla f (x) d(x)^\alpha w(x) \|_{L^1(\O)}
$$

In the case $\alpha=1$, bound (\ref{intfrac}), as before, by

$$
C \int_\O Mg(x) d(x) |\nabla f(x)| \, dx \le C \| Mg(x) w^{-1}(x)\|_{L^{\infty}(\O)} \| \nabla f(x) d(x) w(x)\|_{L^1(\O)}
$$
and the result follows, since by \cite[Theorem 4]{Mu}, if $w\in A_1$,
$$
\| Mg(x) w^{-1}(x)\|_{L^{\infty}(\O)} \le C \| g(x) w^{-1}(x)\|_{L^{\infty}(\O)}
$$

\Box

\begin{remark}
If a weight $w$ satisfies condition (\ref{cond-1peso-L1}), then
$w^q$ belongs to the class $A_1$.
\end{remark}

\section{The two-weighted case}

\begin{theorem}
\label{teo41}

Let $\O\subseteq\mathbb{R}^n$ be a bounded John domain. If $w_1$
and $w_2$  are nonnegative functions such that there exists a
constant $K<\infty$ such that

\begin{equation}
\label{cond-2pesos} |Q|^{\frac1n -1}\left(\int_Q w_1(x) \, dx
\right)^{1/q} \left( \int_Q w_2(x)^{1-p'} \, dx \right)^{1/p'} \le
K
\end{equation}
and $w_1, w_2^{1-p'}$ satisfy the following `reverse
doubling' condition:

\begin{equation}
\mbox{ for any } \epsilon \in (0,1)\mbox{ there exists } \delta\in (0,1) \mbox{ such that }
\int_{\epsilon Q} w(x) \, dx \le \delta \int_Q w(x) \, dx \qquad
\end{equation}
where $Q$ is any $n$-dimensional cube, and $K$ is independent of
$Q$, then

$$
\inf_{a\in \mathbb{R}}\| (f(x)-a) w_1^{1/q}(x)\|_{L^q(\Omega)} \le C
\|\nabla f(x) d(x)^\alpha w_2(x)^{1/p}\|_{L^p(\O)}
$$
for all locally Lipschitz $f$, whenever
$1< p<q<\infty$ and $\alpha \in [0,1]$. If $p=q$, condition (\ref{cond-2pesos}), should be
replaced by requiring that there exist $r>1$ such that

\begin{equation}
\label{cond-2pesosb} |Q|^{\frac{\alpha}{n}+ \frac1q
-\frac1p}\left(\frac{1}{|Q|} \int_Q w_1(x)^r \, dx \right)^{1/qr}
\left( \frac{1}{|Q|} \int_Q w_2(x)^{(1-p')r} \, dx \right)^{1/p'r}
\le K(r)
\end{equation}
\end{theorem}

\proof
As in the previous theorems, by duality it suffices to bound
$\int_\O (f-f_\varphi)(y) g(y) \, dy$ for any $g$ such that
$\|g(x) w(x)^{-1/q} \|_{L^{q'}}< \infty$.

We begin by the case $\alpha\in [0,1)$. Using the bound (\ref{intfrac}) and H\"older's inequality, we
obtain

$$
 \int_\O |(f- f_\varphi)(y) g(y)| \, dy \le  C
\| \nabla f(x) d(x)^\alpha w_2(x)^{1/p} \|_{L^p(\O)} \|
I_{1-\alpha} |g(x)| w_2(x)^{-1/p}\|_{L^{p'}(\O)}
$$

But, by condition (\ref{cond-2pesos}) (respectively, condition
(\ref{cond-2pesosb})) and \cite[Theorem 1]{SW},

$$
\|I_{1-\alpha} |g(x)| w_2^{-1/p} \|_{L^{p'}} \le C \|g(x)
w_1(x)^{-1/q} \|_{L^{q'}}
$$
as we wanted to show.

In the case $\alpha=1$, bound (\ref{intfrac}), as before, by

$$
C \int_\O Mg(x) d(x) |\nabla f(x)| \, dx \le C \| Mg(x) w_2^{-1}(x)\|_{L^{p'}(\O)} \| \nabla f(x) d(x) w_2(x)\|_{L^p(\O)}
$$
and the result follows, since by condition (\ref{cond-1peso}) and \cite[Theorem 1.2]{CU},
$$
\| Mg(x) w_2^{-1}(x)\|_{L^{p'}(\O)} \le C \| g(x) w_1^{-1}(x)\|_{L^{q'}(\O)}
$$

\Box

 \begin{remark}
In the previous theorem we may assume that $q\le \frac{np}{n- p(1-\alpha)}$ (and thus $p(1-\alpha)<n$),
since otherwise $w_1$ equals zero almost everywhere on $\{ w_2< \infty \}$.
This was observed in \cite[Remark b]{S2}.
\end{remark}

\begin{theorem}
\label{teo42}

Let $\O\subseteq\mathbb{R}^n$ be a bounded John domain. If $w_1$
and $w_2$  are nonnegative functions such that there exists a
constant $K<\infty$ such that

\begin{equation}
\label{cond-2pesos-L1}
M(w_2(x)) \le w_1(x)
\end{equation}
for almost all $x$, then

$$
\inf_{a\in \mathbb{R}}\| (f(x)-a) w_1(x)\|_{L^1(\Omega)} \le C
\|\nabla f(x) d(x) w_2(x)\|_{L^1(\O)}
$$
for all locally Lipschitz $f$.
\end{theorem}

\proof
By duality, it suffices to bound $\int_\O (f-f_\varphi)(y) g(y)
\, dy$ for any $g$ such that $\|g(x) w_1^{-1}(x)\|_{L^\infty(\O)}<
\infty$.

As before, bound (\ref{intfrac}) by

$$
C \int_\O Mg(x) d(x) |\nabla f(x)| \, dx \le C \| Mg(x) w_2^{-1}(x)\|_{L^\infty (\O)} \| \nabla f(x) d(x) w_2(x)\|_{L^1(\O)}
$$
and the result follows, since by condition (\ref{cond-1peso-L1}) and \cite[Theorem 4]{Mu},
$$
\| Mg(x) w_2^{-1}(x)\|_{L^\infty(\O)} \le C \| g(x) w_1^{-1}(x)\|_{L^\infty(\O)}
$$

\Box

\begin{remark}

Notice that if one wanted to prove the more general inequality

$$
\inf_{a\in \mathbb{R}}\| (f(x)-a) w_1^{\frac{n-1+\alpha}{n}}(x)\|_{L^{\frac{n}{n-1+\alpha}}(\Omega)} \le C
\|\nabla f(x) d(x)^\alpha w_2(x)\|_{L^1(\O)}
$$
following the proof of the one-weighted case, one would need to know that, if $E_t=\{|I_{1-\alpha}g(x)| > t\}$, then

$$
\int_{E_t} w_1(x) \, dx \le C t^{-\frac{n}{n-1+\alpha}} \left( \int |g(x)| w_2(x) \, dx\right)^{\frac{n}{n+\alpha-1}}.
$$

Unfortunately, we were unable to find neither proof of this inequality under the conditions of the previous theorem (or any other sufficient
conditions on the weights $w_1$, $w_2$) nor any counterexample to the required weak inequality.
Such a result is beyond the scope of this paper, but it is worth noticing that it would immediately imply the above two-weighted
Sobolev-Poincar\'e inequality which would complete Theorem \ref{teo41} in the case $p=1$.

\end{remark}

\bigskip

{\bf Acknowledgement.} We wish to thank Jos\'e Sabina de Lis for giving us reference \cite{Jo}.

\end{document}